\documentclass{amsart}

\usepackage{siunitx}
\usepackage{pifont}
\usepackage{latexsym}
\usepackage{graphicx}
\usepackage{multicol,multirow}
\usepackage{amsmath,amssymb,amsfonts}
\usepackage{mathrsfs}
\usepackage{amsthm}
\usepackage{rotating}
\usepackage{appendix}
\usepackage[numbers,sort&compress]{natbib}
\usepackage{ifpdf}
\usepackage[T1]{fontenc}
\usepackage{times}
\usepackage{newtxmath}
\usepackage{textcomp}%
\usepackage{xcolor}%
\usepackage{hyperref}
\usepackage{lipsum}
\usepackage{url}
\usepackage{enumitem}
\usepackage{tablefootnote}

\usepackage{algorithm} 

\usepackage[linesnumbered,ruled,algo2e]{algorithm2e}

\DeclareMathOperator{\diag}{diag}

\DeclareMathOperator*{\argmin}{arg\,min}

\begin{document}

\title[Fast PCA for Cryo-EM Images]{Fast Principal Component Analysis for Cryo-EM Images}

\author[N. F. Marshall]{Nicholas F. Marshall}
\email{marsnich@oregonstate.edu}
\author[O. Mickelin]{Oscar Mickelin}
\email{hm6655@princeton.edu}
\author[Y. Shi]{Yunpeng Shi} 
\email{yunpengs@princeton.edu}
\author[A. Singer]{Amit Singer}
\email{amits@math.princeton.edu}

\thanks{
A.S. was supported by grants from the AFOSR FA9550-20-1-0266; the Simons Foundation Math+X Investigator Award; NSF BIGDATA Award IIS-1837992; NSF DMS-2009753; and NIH/NIGMS
1R01GM136780-01. We are grateful to Kenny Huang for running some numerical experiments in the initial stage of this project. N.F.M thanks the IPAM 
long program on Computational Microscopy
for hosting him 
during the preparation of this manuscript.}

\newcommand{\cmark}{{\color{green}\ding{51}}}%
\newcommand{\xmark}{{\color{red}\ding{55}}}%


\begin{abstract}
Principal component analysis (PCA) plays an important role in the analysis of cryo-EM images for various tasks such as classification, denoising, compression, and ab-initio modeling. We introduce a fast method for estimating a compressed representation of the 2-D covariance matrix of noisy cryo-electron microscopy projection images  that enables fast PCA computation. Our method is based on a new algorithm for expanding images in
the Fourier-Bessel basis (the harmonics on the disk), which provides a convenient way to handle the effect of the contrast transfer functions. For $N$ images of size $L\times L$, our method has time complexity $O(N L^3 + L^4)$ and space complexity $O(NL^2 + L^3)$. In contrast to previous work, these complexities are independent of the number of different contrast transfer functions of the images. We demonstrate our approach on synthetic and experimental data and show acceleration by factors of up to two orders of magnitude.
\end{abstract}

\maketitle

\section{Introduction}

We study the problem of computing  a compressed representation of the covariance matrix of 2-D cryo-EM images for the purpose of performing principal component analysis (PCA).
More precisely, 
we consider an image formation model where the measurement $g_i$ is defined by
\begin{equation} \label{model}
g_i = h_i * f_i + \varepsilon_i \quad \text{for} \quad i=1,\ldots,N,
\end{equation}
where $h_i$ is a radial function, $*$ denotes convolution, $f_i$ is an image, and  $\varepsilon_i$ the noise term. 

We are motivated by single particle cryo-electron microscopy (cryo-EM) imaging, which is an important technique for determining the 3-D structure of macromolecules. In particular, the single particle reconstruction (SPR) problem asks to recover the 3-D structure of a macromolecule from noisy 2-D images of its tomographic projections along unknown viewing angles. In cryo-EM, the mathematical model is a special case of \eqref{model} and is of the form
\begin{eqnarray}\label{EQN:full_defn}
    g_i(x') = h_i \ast \int_\mathbb{R} \varphi_i \left(R_i^{-1} x \right) dx_3 + \varepsilon_i(x'), \qquad i = 1, \ldots , N,
\end{eqnarray}
where $x = (x',x_3) \in \mathbb{R}^2 \times \mathbb{R} \cong \mathbb{R}^3$ are 3-D spatial coordinates,  $h_i$ is the point spread function, $\varphi_i : \mathbb{R}^3 \rightarrow \mathbb{R}$ is the electrostatic potential of a molecule, $R_i \in SO(3)$ is a 3-D 
rotation, and $\varepsilon_i$ the noise term. 
In computational microscopy, it is typical to work with the Fourier transform of the point spread function, which is known as the contrast transfer function (CTF).
In the simplest case, each measurement could correspond to a single fixed molecule potential function  $\varphi_i = \varphi$; however, in general, we may assume that each $\varphi_i$ could be a random variable representing a mixture of molecules, conformational heterogeneity, cases where the images are not perfectly centered, or other measurement imperfections \cite{liu1995estimation,penczek2011identifying}.

In general, each measurement $g_i$ can be associated with a different point spread function; however, in practice,  a group of measurements, called a defocus group, can share a common point spread function.
We assume that the measurements are grouped into $M \le N$ defocus groups. Given $g_i$ and $h_i$ for $i = 1,\ldots,N$, our goal is to estimate the 2-D covariance function $c : \mathbb{R}^2 \times \mathbb{R}^2 \rightarrow \mathbb{R}$ of the images
\begin{equation} \label{covariance}
c(x',y') := \mathbb{E} \left[ \left(f(x') - \bar{f}(x') \right) \left( f(y') - \bar{f}(y') \right)\right],
\end{equation}
where $f$ is a random variable from the same distribution as the images $f_i$, and
$\bar{f}(x') = \mathbb{E}[ f(x')]$. We assume that the distribution of the images is invariant to in-plane rotations (which is typically the case in cryo-EM).

 In cryo-EM,  the random variable $f$ is of the form $f(x') = \int_{\mathbb{R}} \Phi(R^{-1} x) dx_3$, where the random variable $R$ is an unknown viewing angle, and the random variable $\Phi$ is a molecule potential. There is generally no physical reason for a molecule to prefer one in-plane rotation to another so distributions of random variables of this form are generally invariant to in-plane rotations. In the field of cryo-EM processing, the covariance function $c$ is simply referred to as the 2-D covariance.

\subsection{Motivation}
The 2D-covariance is an essential component of a number of computational techniques in cryo-EM; we survey a few of these below.

First, we are motivated by PCA, which is a ubiquitous technique in statistics, data science, and computational mathematics and has applications to dimensionality reduction, denoising, visualization, among others. The principal components (that is, the top eigenvectors of the digitized covariance matrix) have a number of uses in the computational cryo-EM pipeline. The subspace corresponding to the top eigenvectors of the covariance matrix identifies salient features of the dataset which enables, for instance, improved methods for image classification and visualization, such as Multivariate Statistical Analysis \cite{vanHeel,van1984multivariate,van2009multivariate,van2016multivariate}. These techniques improve computational speed, since clustering becomes computationally easier in a space of reduced dimension, as well as accuracy, since dimensionality reduction by PCA amplifies the effective signal-to-noise ratio because many coordinates for which noise dominates the signal are eliminated \cite[Table~3]{zhao2014rotationally}.

Second, the covariance matrix has applications in the method of moments, a classical statistical inference method, applied to cryo-EM \cite{kam1980reconstruction}. In this method, the 2-D covariance is used to compute the similarly defined autocorrelation function of the underlying 3-D structure. Under further assumptions such as sufficient non-uniformity of the distribution of the viewing angles \cite{sharon2020method} or sufficient sparsity of the molecular structure \cite{bendory2022autocorrelation}, this autocorrelation function determines the 3-D density map either up to a finite list of possible structures or uniquely, respectively. This has been further developed into principled methods for ab initio estimation of cryo-EM structures \cite{sharon2020method,bendory2022autocorrelation}, with reduced risk of user-induced model bias in the initial model. Alternatively, when additional information is available, for instance one \cite{huang2022orthogonal} or two \cite{levin20183d} noiseless projection images, or the 3-D structure of a related, homologous structure \cite{bhamre2015orthogonal,bhamre2017anisotropic}, the 3-D density map is uniquely determined by the autocorrelation, without requiring any structural assumptions.

Third, the covariance matrix has applications to denoising and CTF-correcting projection images. Covariance Wiener Filtering (CWF) \cite{CWF} is an approach which uses the classic Wiener filtering framework with the estimated covariance matrix to solve the image deconvolution and denoising problem. The technique represents images in a lower dimensional subspace that is formed from PCA using the estimated covariance matrix. The method then applies Wiener filtering to correct the CTFs and denoise the images in this reduced subspace.

Compared to the standard PCA problem, the cryo-EM setting exhibits further computational challenges, since the estimation method also has to account for 
convolution with the point spread function, which destroys information of the resulting convolved function; see \S\ref{sec:technical_details} for a more precise statement. On the other hand, the problem has additional symmetries making fast algorithms possible.
In this paper, we present a new fast algorithm for estimating the covariance matrix that improves upon past approaches (especially when there are a large number of defocus groups) in terms of time and space complexity.

\subsection{Main contribution}
The main contribution of this paper is a new computational method for estimating the covariance \eqref{covariance} from $N$ measurements of the form \eqref{model} encoded by $L \times L$ digitized images. The presented fast method has time complexity $O(N L^3 + L^4)$ independent of the number $M \le N$ of defocus groups. This is in contrast to past methods, where this complexity scales poorly with $M$ and involves $O(M T L^4+N L^3)$ operations \cite{CWF}, where $T$ is the number of iterations needed in a Conjugate Gradient step. Many modern cryo-EM experimental datasets fall into the computationally challenging regime where $M$ scales with $N$.

Our fast method hinges on a new fast and accurate method for expanding $L \times L$ images into the Fourier-Bessel basis, which provides a convenient way to handle convolution of radial functions (such as point spread functions) with images:
namely, convolution with radial functions can be expanded as a diagonal operator operating on the basis coefficients \cite{FLE}.

The Fourier-Bessel basis functions are \emph{harmonics} on the disk: the standing waves associated with the 
resonant frequencies of a disk shaped drum with a fixed boundary. More precisely, the harmonics on the disk are eigenfunctions of the Laplacian on the unit disk that satisfy Dirichlet boundary conditions. In computational mathematics, this basis is referred to as the Fourier-Bessel basis, since the basis functions can be expressed as a product of a Bessel function and a complex exponential; see \eqref{eigenfun} for a definition.

Because of this simple structure, the covariance matrix of clean images can be estimated by a simple closed-form solution, without using the (computationally expensive) conjugate gradient method from previous approaches. Simultaneously, the covariance matrix retains its block diagonal structure, meaning that its diagonal blocks can be estimated separately and independently, which altogether makes PCA fast.

We present numerical results of covariance estimation on 
 synthetic and experimental data. Additionally, we show how the estimated covariance matrix can be used to denoise images using CWF, and perform PCA to visualize eigenimages from experimental data. Code implementing the method is publicly available online\footnote{Code is available at \url{https://github.com/yunpeng-shi/fast-cryoEM-PCA}}. Moreover, our approach has the potential to generalize to settings beyond cryo-EM, where PCA is used for signals estimated under more general group actions \cite{bandeira2017estimation}. 

\subsection{Organization}
The remainder of the article is organized as follows. In \S\ref{sec:methodology}, we describe the computational method. In \S \ref{sec:comp_synth}, we present numerical results for synthetic data. In \S \ref
{sec:comp_real} we present numerical results for experimental data. In \S \ref{discuss}, we discuss the results and possible extensions.

\section{Methodology}\label{sec:methodology}

\subsection{Notation}
For two $M\times N$-matrices $A$ and $B$, we denote their Hadamard (or entrywise) product by $A\odot B$, the Hadamard division of $A$ and $B$ by $A \oslash B$ and the $\ell$th Hadamard power of $A$ by $A^{\odot \ell}$. These operations are defined elementwise by
\begin{equation}
    (A\odot B)_{jk} = A_{jk}B_{jk}, \quad  (A\oslash B)_{jk} = \frac{A_{jk}}{B_{jk}}, \quad (A^{\odot \ell})_{jk} = A_{jk}^\ell, 
\end{equation}
respectively.  If $w$ is an $N$-dimensional vector, then $\diag(w)$ denotes the $N \times N$ matrix with $w$ along its diagonal, i.e., $\diag(w)_{jj} = w_j$, and zeros elsewhere. If $f : \mathbb{R}^2 \rightarrow \mathbb{R}$ is a radial function, we write $f(x) = f(|x|)$ to mean that $f$  can be expressed as a function only of the magnitude $|x|$ of $x$.

\subsection{Technical details}\label{sec:technical_details}

We make the following assumptions 
\begin{enumerate}[label=(A\arabic*)]
\item \label{A1}
We assume that the point spread functions $h_i$ in the model \eqref{model} are radial functions; this implies that their Fourier transform (the CTFs) are also radial. In systems where astigmatism is present and the point spread function deviates slightly from a radial function, our approach can be used as an initial approximation that could be refined using the Conjugate Gradient method. 
\item \label{A2} We assume that the underlying images $f_i$ in the model 
\eqref{model} are i.i.d. random variables whose distribution is
invariant to in-plane rotations.
\item \label{A3} We assume a technical condition on the Fourier transform of the point spread functions $h_i$ in the model \eqref{model}. Namely, that
the Fourier transforms $\hat{h}_i$ of the $h_i$ satisfy
$$
\inf_{ |\xi|,|\eta| \in [\lambda_\text{min}, \lambda_\text{max}]} \sum_{i=1}^N |\widehat{h}_i(\xi)|^2 |\widehat{h}_i(\eta)|^2  \ge \delta,
$$
where $\delta > 0$ is fixed and $[\lambda_\text{min},\lambda_\text{max}]$ is the interval of Fourier space used in the disk harmonic expansion, see \cite[\S2.4]{FLE}.
\end{enumerate}

\begin{figure}[t]
\centering
\includegraphics[width=0.55\textwidth]{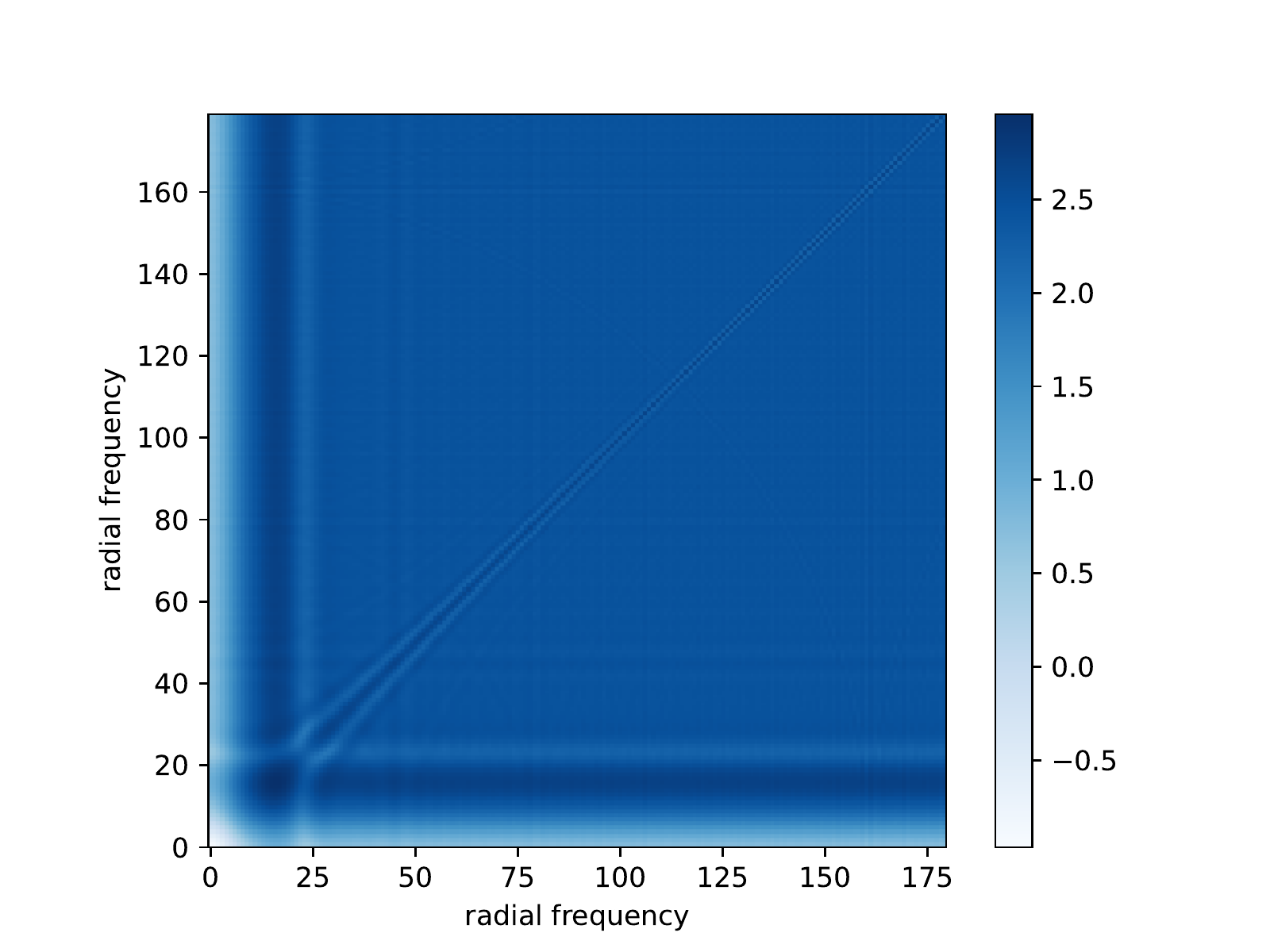}
\caption{We visualize $\log_{10} \left( \sum_{i=1}^N |\widehat{h}_i(|\xi|)|^2 |\widehat{h}_i(|\eta|)|^2\right)$ for each pair of radial frequencies $|\xi|, |\eta|$ for the experimental dataset EMPIAR-10028 \cite{10028} obtained from the Electron Microscopy Public Image Archive \cite{empiar}. All values are greater than $-1$ in the log scale}\label{fig:assumption}
\end{figure}

Informally speaking, assumption \ref{A3} can be interpreted as saying that for any pair of frequencies  $\xi$ and $\eta$ there
is a point spread function $h_i$ whose Fourier transform $\widehat{h}_i$ does not vanish  at $\xi$ or $\eta$.
 Assumption \ref{A1} implies that $\widehat{h}_i(\xi) = \widehat{h}_i(|\xi|)$ is a radial function.
Figure \ref{fig:assumption} shows the values of $\sum_{i=1}^N |\widehat{h}_i(\xi)|^2 |\widehat{h}_i(\eta)|^2$ for each pair of radial frequency $\xi, \eta$ in log scale, for 1081 distinct CTF images of size $L=360$, whose defocus values range from 0.81 to 3.87 \si{\mu m}, for an experimental dataset; see \S\ref{sec:comp_real} for more details.

This assumption is much weaker than assuming that for all $i$ that $|\widehat{h}_i|$ does not vanish at any frequency. In the latter case, we could just use $h_i$ to invert each equation to get access to the underlying functions $f_i$.
If we had direct access to the underlying functions $f_i$, then we could approximate the covariance function $c$ by the sample covariance 
$$
c_N(x',y') := \frac{1}{N-1} \sum_{i=1}^N [f_i(x') - \tilde{f}(x')] \cdot [f_i(y') - \tilde{f}(y')],
$$
where $\tilde{f}(x') = \sum_{i=1}^N \frac{1}{N} f_i(x')$ is the sample mean.
Indeed, $c_N(x',y') \rightarrow c(x',y')$ by the law of large numbers. To clarify why it is useful for $\widehat{h}_i$ not to vanish, 
note that in the Fourier domain our measurement model can be expressed as
$$
\widehat{g}_i = \widehat{h}_i \widehat{f}_i + \widehat{\varepsilon}_i,
$$
where $\widehat{g}_i,\widehat{h}_i,\widehat{f}_i$, and $\widehat{\varepsilon}_i$ denote the Fourier transforms of $g_i,h_i,f_i,$ and $\varepsilon_i$, respectively. Since taking the Fourier transform changes convolution to point-wise multiplication, if each Fourier transform $|\widehat{h}_i| \ge \delta$ for some $\delta > 0$, then we could estimate $c$ by first estimating $\widehat{f}_i$ from $\widehat{g}_i$ for $i =1,\ldots,N$, and then using the sample covariance matrix. However,  in practice, the CTF $\widehat{h}_i$ is approximately a radial function with many zero-crossings, which means that multiplication by $\widehat{h}_i$ destroys information in the corresponding frequencies, making the restoration from a single image ill-posed.

\subsection{ Fourier-Bessel}\label{sec:def_FB}
The main ingredient in our fast covariance estimation is a fast transform into a convenient and computationally advantageous basis, known as the Fourier-Bessel basis (which consists of the harmonics on the disk: eigenfunctions of the Laplacian on the disk that obey Dirichlet boundary conditions). This specific choice of basis has a number of beneficial properties:
\begin{enumerate}[label=(\roman*)]
\item it is orthonormal, 
\item it is ordered by frequency,
\item it is steerable i.e., images can be rotated by applying a diagonal transform to the basis coefficients,
\item it is easy to convolve with radial functions i.e., images can be convolved with radial functions by applying a diagonal transform to the basis  coefficients.
\end{enumerate}
These properties have made the Fourier-Bessel basis a natural choice in a number of imaging applications \cite{shi2022ab,CWF,rangan2020factorization,zhao2014rotationally,zhao2016fast} and will be central to the development of our fast covariance estimation method. In polar coordinates $(r,\theta)$ in the unit disk $\{ (r,\theta) : r \in [0,1), \theta \in [0, 2\pi)  \}$ the Fourier-Bessel basis functions are defined by
\begin{equation} \label{eigenfun}
    \psi_{nk}(r,\theta) = \gamma_{nk} J_n(\lambda_{nk}r)e^{\imath n \theta},
\end{equation}
where $\gamma_{nk}$ is a normalization constant, $J_n$ is the $n$-th order Bessel function of the first kind \cite[\S 10.2]{dlmf}, and $\lambda_{nk}$ is its $k$-th smallest positive root; the indices $(n,k)$ run over $\mathbb{Z} \times \mathbb{Z}_{> 0}$.

Recent work \cite{FLE} has devised a new fast algorithm to expand $L\times L$-images into $ \sim L^2$ Fourier-Bessel basis functions.
Informally speaking, given $ \sim L^2$ basis coefficients, the algorithm can evaluate the function on an $L \times L$ grid in $O(L^2 \log L)$ operations; the adjoint can be computed in the same number of operations, which makes iterative methods fast.

Compared to previous expansion methods \cite{zhao2014rotationally,zhao2016fast}, it enjoys both theoretically guaranteed accuracy and lower time complexity. 
\subsection{Key property of Fourier-Bessel basis}\label{sec:key}

A key property of the Fourier-Bessel basis is that convolution with radial functions are diagonal transformations in any truncated basis expansion.  More precisely, the following result holds \cite[Lemma 2.3]{FLE}: suppose that $f = \sum_{(n,k) \in I} \alpha_{nk} \psi_{nk}$ for some index set $I$, and $h = h(|x|)$ is a radial function. Then,
\begin{equation} \label{diagt}
P_\mathcal{I} (f*h) = \sum_{(n,k) \in I} \alpha_{nk}\widehat{h}(\lambda_{nk}) \psi_{nk}
\end{equation}
where $P_\mathcal{I}$ denotes orthogonal projection onto the span of $\{\psi_{nk}\}_{(n,k)\in I}$, $\widehat{h}$ is the Fourier transform of $h$, and $\lambda_{nk}$ is the $k$-th positive root of $J_n$.

We emphasize that the weights of the diagonal transform in \eqref{diagt} are not the coefficients of $h$ in the disk harmonic expansion. Indeed, 
since $h$ is radial, it follows from \eqref{eigenfun} that the coefficients $\beta_{nk}$ of $h$ in the basis $\psi_{nk}$ satisfy $\beta_{nk} = 0$ when $n \not =0$. Computing the weights $\widehat{h}(\lambda_{nk})$ from the coefficients $\beta_{nk}$ of $h$ in the basis, would involve computing weighted sums of the Fourier transforms of the basis functions: $\widehat{h}(\lambda_{nk}) = \sum_{l \in \mathcal{J}} \beta_{0l} \widehat{\psi}_{0l}(\lambda_{nk})$, for some index set $\mathcal{J}$.

As an alternative to the disk harmonic basis expansion, one can consider simply taking the Fourier-transform of \eqref{model}, which also leads to a diagonal representation of the convolution operator. However, the discrete Fourier transform does not have the steerability property, which is essential for the covariance estimation. Another attempt could be to use the polar Fourier transform. However, this representation is not invariant to arbitrary in-plane rotations, but only to finitely many rotations as determined by the discretization spacing of the grid. These expansions are therefore unsuitable for the goal of this article and we instead use expansions into the Fourier-Bessel basis, although other steerable bases could be considered \cite{landa2017approximation,landa2017steerable}. Table~\ref{tab:basis_properties} summarizes the considerations that make the Fourier-Bessel basis a natural choice of basis.

\begin{table}
    \caption{Summary of desirable properties of a few different basis candidates.}
    \label{tab:basis_properties}
    \centering
    \begin{tabular}{c|cccc}
       Basis  & Orth. & Cont. steerable &  Radial convolution
       diag. & Fast expansion\tablefootnote{The basis expansion from Cartesian grid representation can be completed within $O(L^2)$ operations up to log factors.} 
       \\
      \hline
    Real &\cmark & \xmark &\xmark & \cmark\\  
    2-D discrete Fourier &\cmark & \xmark &\cmark &\cmark\\  
        Polar Fourier &\xmark & \xmark  \, (discrete) &\cmark & \cmark\\  
        
          PSWF \cite{landa2017approximation,landa2017steerable} &\cmark & \cmark &\xmark &\xmark\\      
      Fourier-Bessel \cite{FLE}\tablefootnote{The new expansion algorithm \cite{FLE} improves the previous computational method \cite{zhao2014rotationally} in terms of accuracy guarantees, computational complexity and the fact that it derives weights such that radial convolution is a diagonal operation.} & \cmark & \cmark &\cmark & \cmark\\
      \hline
    \end{tabular}
\end{table}

\subsection{Block diagonal structure}
The steerable property of the Fourier-Bessel basis implies that the 2D-covariance will be block diagonal in this basis.
 A full representation of an $L^2\times L^2$ matrix requires $O(L^4)$ elements, but this is reduced to $O(L^3)$ nonzero entries by the block diagonal structure. This block diagonal structure follows from the form of the basis functions and that the distribution of in-plane rotations is assumed to be uniform. Indeed, suppose that $f = \sum_{(n,k) \in I} \alpha_{nk} \psi_{nk}.$
For simplicity assume that $f$ has mean zero; subtracting the mean will only change the radial components, since the other components corresponding to non-vanishing angular frequencies have zero mean by merely averaging over all possible in-plane rotations. By assumption~\ref{A2}, the covariance function in polar coordinates satisfies
\begin{equation}\label{eq:cov_in_plane}
    c( (r,\theta), (r',\theta')) = c( (r,\theta + \varphi), (r',\theta' + \varphi)),
\end{equation}
for all $\varphi$ in $[0,2\pi]$. The covariance function in~\eqref{covariance} can be expanded in a double Fourier-Bessel basis expansion as
\begin{equation}\label{eq:c_fle_block_diag}
   c\left((r,\theta),(r',\theta')\right) = \sum_{(n,k) \in I} \sum_{(n',k') \in I} C_{(n k, n' k')}  \psi_{nk}(r,\theta) \overline{\psi_{n'k'}}(r',\theta'),
\end{equation}
where $C_{(n k, n' k')}$ is the covariance matrix in the Fourier-Bessel basis. Combining \eqref{eq:cov_in_plane} and  \eqref{eq:c_fle_block_diag} and integrating $\varphi$ over $[0,2\pi]$ gives

\begin{equation} \label{eq:calcorth}
\begin{split}
    c (r,\theta), (r',\theta')) &= \frac{1}{2\pi}\int_{0}^{2\pi} c( (r,\theta + \varphi), (r',\theta' + \varphi)) d\varphi \\
    &= \sum_{(n,k) \in I} \sum_{(n',k') \in I} C_{(n k, n' k')}
    \frac{1}{2\pi}
    \int_0^{2\pi}\psi_{nk}(r,\theta+\varphi) \overline{\psi_{n'k'}}(r',\theta'+\varphi) d\varphi \\
    &= \sum_{(n,k) \in I} \sum_{(n',k') \in I} C_{(n k, n' k')}  \psi_{nk}(r,\theta) \overline{\psi_{n'k'}}(r',\theta')
    \frac{1}{2\pi}
    \int_0^{2\pi} e^{\imath (n-n')\varphi} d\varphi, \\
    &= \sum_{(n,k) \in I} \sum_{(n',k') \in I} \delta_{n=n'} C_{(n k, n' k')}  \psi_{nk}(r,\theta) \overline{\psi_{n'k'}}(r',\theta'),
    \end{split}
\end{equation}
where $\delta_{n=n'}$ is a Dirac function that is equal to $1$ if $n=n'$ and zero otherwise, and  we note that the 
second to last equality uses the fact that $\psi_{nk}(r,\theta) = \gamma_{nk}J_n(\lambda_{nk}r)e^{\imath n \theta}$.  
Since the coefficients $C_{(nk,n'k')}$ in the expansion \eqref{eq:c_fle_block_diag} are unique, it follows from \eqref{eq:calcorth} that $C_{(nk,n'k')} = 0$ when $n \not = n'$. Hence, the covariance matrix $C_{(n k,n'k')}$ has a block diagonal structure whose blocks consist of the indices $(nk,nk')$ for a given value of $n$.
In the following section, we show how these properties enable a fast method to estimate the covariance matrix.

\subsection{Covariance Estimation}\label{sec:sec1_covar}
In the Fourier-Bessel basis, \eqref{model} is written as
\begin{equation}\label{eq:def_fb_conv}
    G_i = H_i \odot F_i + E_i,
\end{equation}
where $G_i$, $F_i$, and $E_i$ are coefficient vectors of $g_i, f_i, \varepsilon_i$ in the Fourier-Bessel basis, respectively and $H_i$ is the vector encoding the convolution operator of \S\ref{sec:key}, i.e., with components $\widehat{h_i}(\lambda_{nk})$. The vectors are $b$-dimensional column vectors, where $b = O(L^2)$ is the number of basis coefficients. We use this simple structure to obtain a closed form expression for the sample covariance matrix of the $F_i$. We estimate this matrix by minimizing the discrepancy between the sample covariance and the population covariance; more precisely, the estimated covariance matrix $\tilde{C}$ is computed by solving the least squares-problem
\begin{equation}\label{eq:lst_sq}
    \tilde{C} = \argmin_C \sum_{i=1}^N \left\| \left( G_i - H_i \odot \tilde{\mu} \right)\left( G_i - H_i \odot \tilde{\mu} \right)^T - \left(C \odot (H_i H_i^T) + \sigma^2 I \right) \right\|^2_F.
\end{equation}
where $\tilde{\mu}$ is defined by
\begin{equation}
    \tilde{\mu} = \argmin_{\mu} \sum_{i=1}^N \| G_i - H_i \odot \mu \|^2,
\end{equation}
whose solution is
\begin{equation}\label{eq:sample_mu}
\tilde{\mu} =
\left( \sum_{i=1}^N H_i \odot G_i \right) \oslash  \left( \sum_{i=1}^N H_i^{\odot 2} \right).
\end{equation}
The least squares-solution of~\eqref{eq:lst_sq} can be determined by the following system of linear equations
\begin{equation}\label{eq:hadamard_division}
   \left( \sum_{i=1}^N H_i^{\odot 2} (H_i^{\odot 2})^T\right) \odot \tilde{C} = \sum_{i=1}^N (H_i H_i^T) \odot B_i - \sigma^2 \sum_{i=1}^N \diag(H_i^{\odot 2})
\end{equation}
where
\begin{equation}\label{eq:def_Ci}
B_i = (G_i - H_i \tilde{\mu})(G_i - H_i\tilde{\mu})^T .
\end{equation}
  It follows that
\begin{equation}\label{eq:cov}
        \tilde{C} = \left( \sum_{i=1}^N \left[ B_i\odot (H_i H_i^T) - \sigma^2 \text{diag}(H_i^{\odot 2})\right] \right) \oslash \left(\sum_{i=1}^N H_i^{\odot 2} (H_i^{\odot 2})^T\right).
\end{equation}
As discussed in \S\ref{sec:key}, the covariance matrix is block diagonal in the  Fourier-Bessel basis. More precisely, the only non-zero elements $\tilde{C}(nk,n'k')$ of the matrix $\tilde{C}$ are those with $n=n'$. Therefore, the matrices in \eqref{eq:def_Ci} and \eqref{eq:cov} need only be calculated for this subset of indices. Since there is a total of $O(L^3)$ of these indices, this reduces the computational complexity compared to computing with the full matrices.

Note that the covariance matrix estimated from \eqref{eq:cov} may not be positive semidefinite due to subtraction of the term $\sigma^2\text{diag}(H_i^{\odot 2})$. Therefore, when running the method in practice it is beneficial to use an eigenvalue shrinkage method. The computational cost of eigenvalue shrinkage for a matrix with our block structure is $O(L^4)$, see for example \cite{donoho2018optimal,CWF}.
For completeness, we include this computational cost of eigenvalue shrinkage in our overall computational complexity.
Informally speaking, the idea of eigenvalue shrinkage
is to replace the term  $\sum_{i=1}^N \left[ B_i\odot (H_i H_i^T) - \sigma^2 \text{diag}(H_i^{\odot 2})\right]$ in \eqref{eq:cov} by $\sum_{i=1}^N \left[ B_i\odot (H_i H_i^T) \right]$, and then shrink and truncate the eigenvalues in a systematic way before Hadamard division, see \cite{donoho2018optimal} for details. The steps of the algorithm are summarized in Algorithm~\ref{alg:cov}.

\begin{algorithm}
\raggedright
    \SetKwInOut{Input}{Input}
    \SetKwInOut{Output}{Output}

    \Input{Observed images $g_i$, radial functions $h_i$, $i = 1, \ldots , N$}
    \Output{Estimated covariance matrix $\tilde{C}$ in the domain of the Fourier-Bessel basis}
    Expand observed images $g_i$ into the Fourier-Bessel basis, with resulting coefficient vector $G_i$ \label{step1}
    
    Compute the vectors $H_i$ with components $\widehat{h_i}(\lambda_{nk})$, representing the action of the CTFs in the Fourier-Bessel basis \label{step2}
    
    Compute sample mean $\tilde{\mu}$ from \eqref{eq:sample_mu} \label{step3}
    
    Use \eqref{eq:def_Ci} to compute the elements $B_i(nk,nk')$ for all $n,k,k'$ \label{step4}
    
    Use \eqref{eq:cov} to compute the elements of the sample covariance matrix $\tilde{C}(nk,nk')$ for all $n,k,k'$, refined using eigenvalue shrinkage \label{step5}
    
    \caption{Fast covariance estimation method}\label{alg:cov}
\end{algorithm}

The complexity of Step~\ref{step1} of the algorithm is $O(NL^2\log L)$. 
The complexity of Step~\ref{step2} is $O(ML^2\log L)$.
The complexity of Step~\ref{step3} is $O(NL^2)$, since the number of basis coefficients is $O(L^2)$. The  complexity of Step~\ref{step4} is $O(NL^3)$, since there are at most $O(L^3)$ non-zero elements with indices $(nk,nk')$,
The complexity of Step~\ref{step5} is $O(NL^3+L^4)$, where the additional term $\mathcal{O}(L^4)$ comes from the computational complexity of eigenvalue shrinkage. Thus, the total complexity of the algorithm is $O(NL^3+L^4)$.

\subsection{Application to image denoising}\label{sec:sec2_wiener}
Given an estimate of the covariance matrix $\tilde{C}$, the CWF approach estimates the $F_i$ by a linear Wiener filter \cite{CWF},
\begin{equation}
    \tilde{F_i} =  \tilde{\mu} + \tilde{C} \diag(H_i)\left(  \tilde{C} \odot (H_iH_i^T) + \sigma^2 I \right)^{-1}(G_i - H_i\odot \tilde{\mu}),
\end{equation}
see \cite{CWF} for more details.

\section{Synthetic data results}\label{sec:comp_synth}

\begin{figure}[t]
\centering
\includegraphics[width=0.8\textwidth]{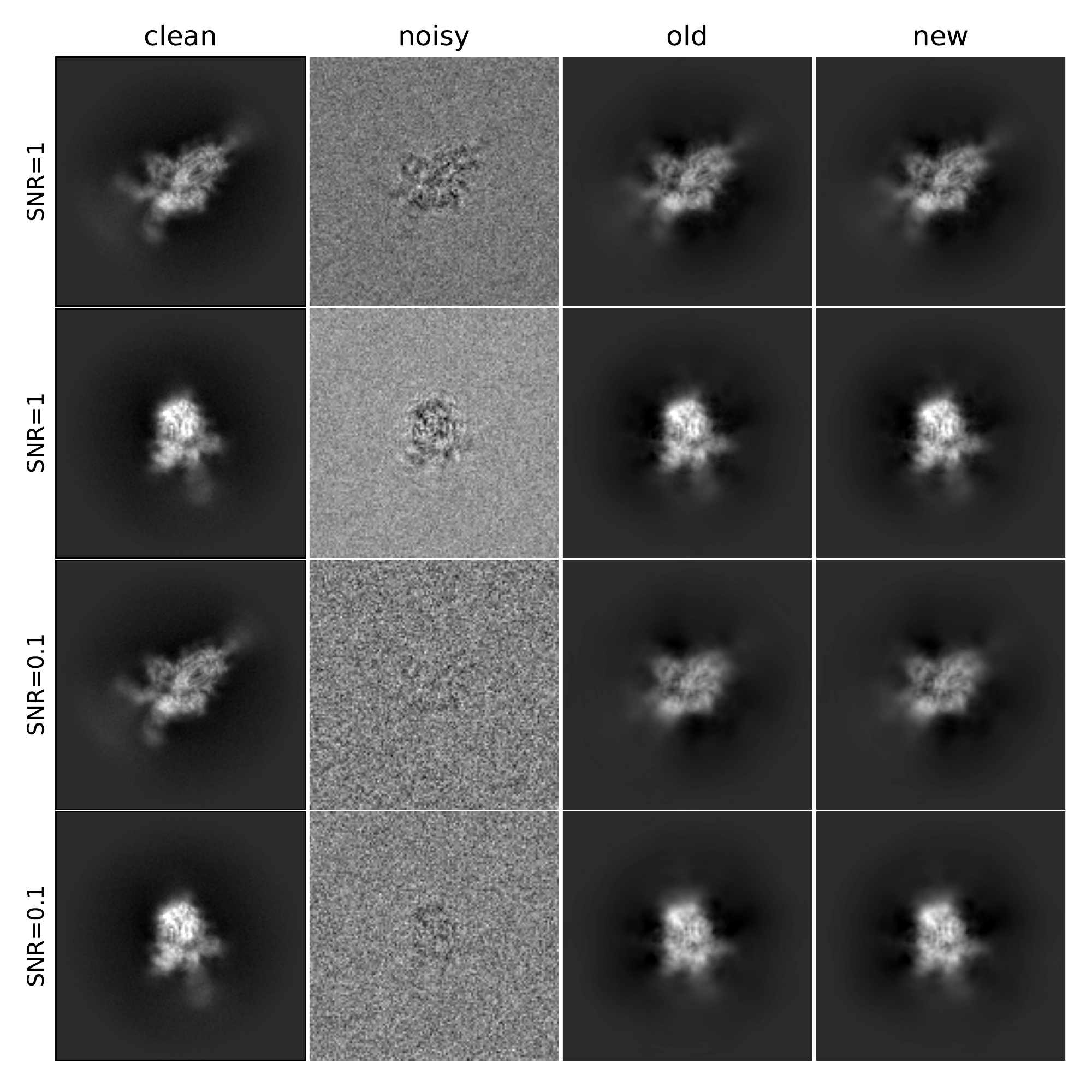}
\caption{Clean, noisy and denoised images. The covariance estimation used $N = 10000$ images, and parameters $L=128$ and $M=100$}\label{fig:synthetic_image}
\end{figure}

We compare the timings of our fast method to previous approaches \cite{CWF}, for synthetic images generated from the 3-D volume of SARS-CoV-2 (Omicron) spike complexes \cite{guo2022structures} (EMD-32743), from the online EM data bank \cite{lawson2016emdatabank}. The original volume has size $512$ in each dimension, with pixel size 0.832 \AA.  We downsample the original volume to size $L\times L \times L$, with $L=32,128$ and 512, respectively, and show the computational times. To generate the synthetic noisy images, we first generate 10000 clean projection images of the 3-D volume from random and uniformly distributed viewing directions. We next divide the set of clean images into a number of defocus groups, where the defocus values range from 1 \si{\mu m} to 4 \si{\mu m}. For all CTFs, we
set the voltage as 300 \si{kV} and the spherical aberration
as 2 \si{mm}. After convolving the images with their CTFs, we add colored noise with power spectral density $1/(r  L/20 + 1)$ up to a constant scale, where $r\in [0,1]$ is the radial frequency. For both the previous method and ours, the CTFs and the noisy images are whitened before estimating the covariance. A few sample clean and noisy images are shown in Figure~\ref{fig:synthetic_image}. All experiments were carried out on a machine with 72 cores. We note that our implementation does not take advantage of all cores since some of the packages we use are not fully optimized for parallel computing; improving this is a technical direction for future work.

\begin{figure}[t]
\centering
\includegraphics[width=.95\textwidth]{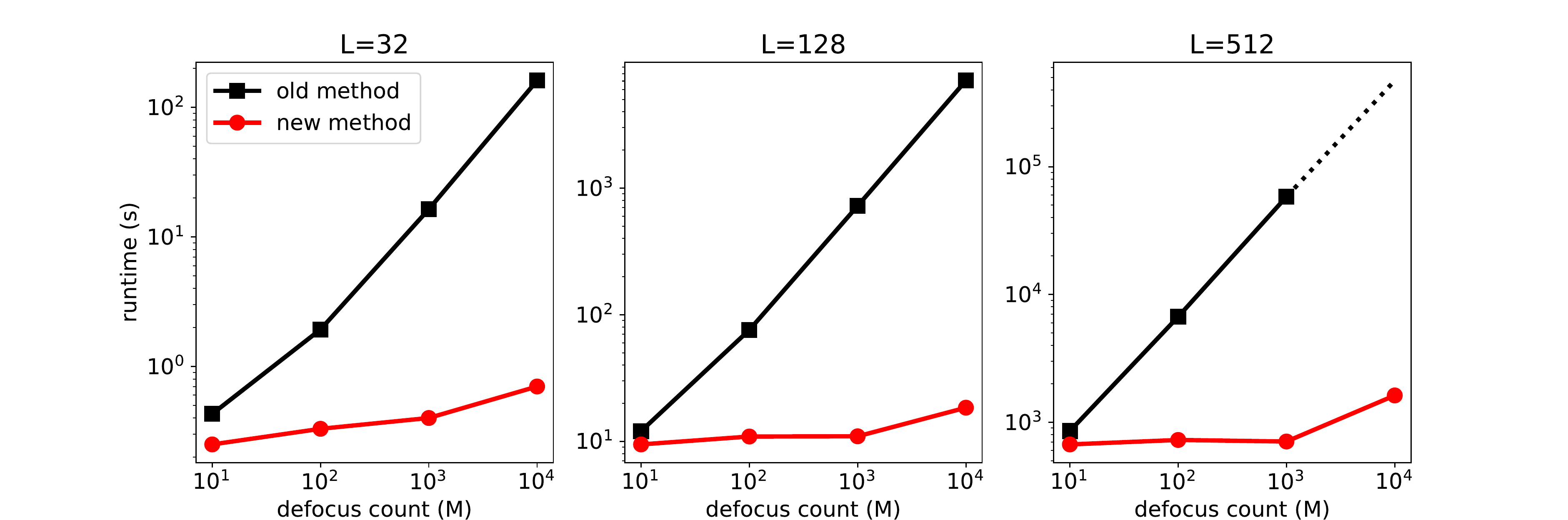}
\caption{Timing comparison for covariance matrix estimation of $10000$ images of size $L \times L$. The old method \cite{CWF} timing for $L=512$ and defocus count $10^4$ is extrapolated due to time and memory constraints}\label{fig:synth_cov}
\end{figure}

Figure~\ref{fig:synth_cov} shows the time required to estimate the covariance matrices as a function of the number of defocus groups. Note that the runtime for the previous approach for the largest image and defocus group sizes is infeasibly large and that our fast method exhibits a speedup of up to three orders of magnitude.

\begin{figure}[htb]
\centering
\includegraphics[width=.95\textwidth]{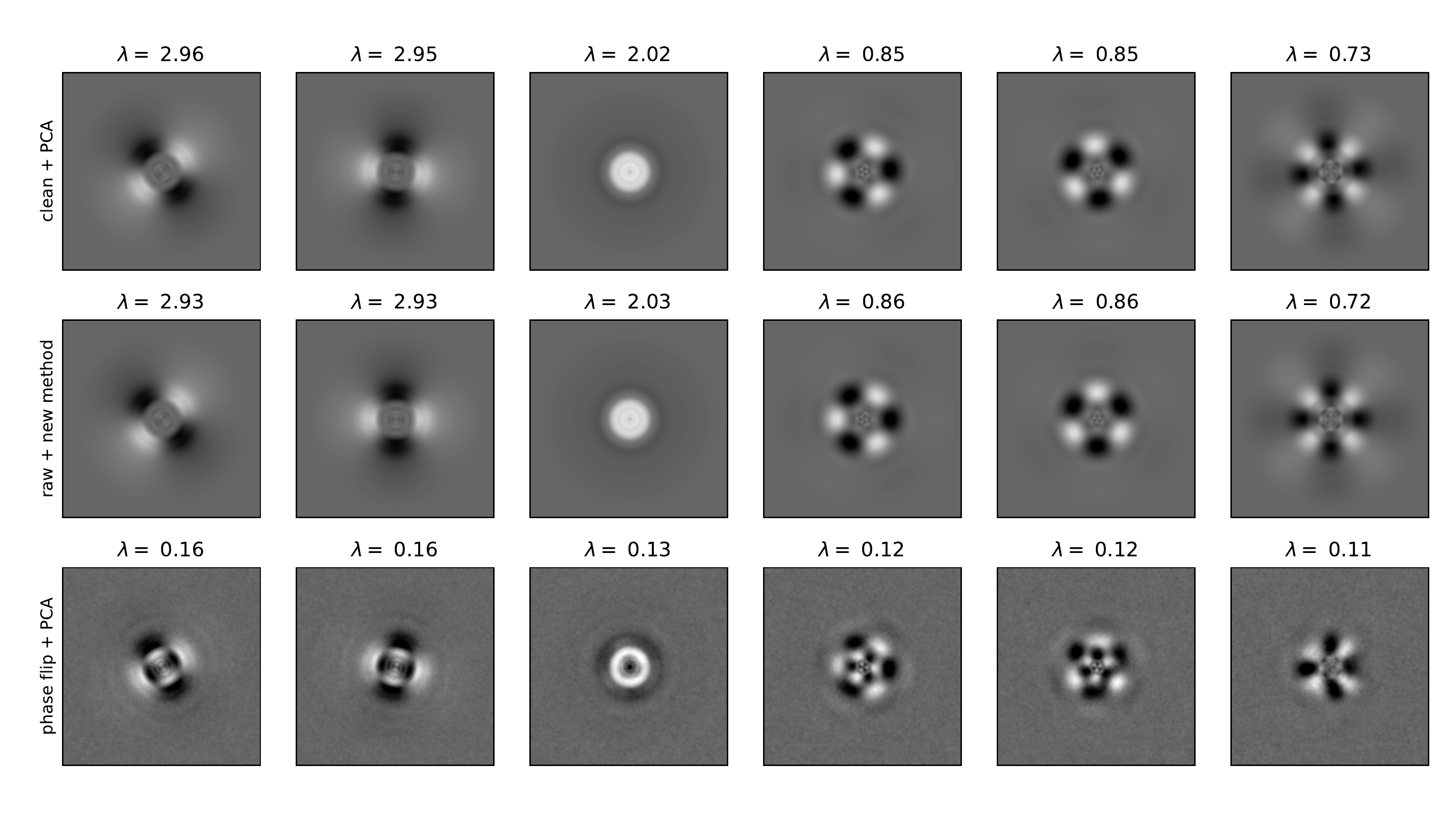}
\caption{Top 6 eigenimages computed by traditional PCA on $10^6$ clean images (top panel), our new method on $10^4$ raw images (middle panel), and traditional PCA on $10^4$ phase-flipped images (bottom panel). The signal-to-noise-ratio for the images for the new method and traditional PCA was 0.1}\label{fig:eig_im}
\end{figure}

Figure~\ref{fig:eig_im} shows the top six principal components estimated by our method and by traditional PCA using $10^4$ raw images, where $L=128$, SNR$=0.1$ and $M=100$ defocus groups, compared to traditional PCA on $10^6$ images. We use the sample covariance matrix of phase-flipped images in real space for the traditional PCA. For all methods, we use $\lambda$ to denote $100$ times the eigenvalues of the eigenimages. The eigenimages from the traditional PCA look much noisier than ours, and contain artifacts that are due to imperfect CTF correction (see e.g. the circular artifacts in top 3 eigenimages in Figure \ref{fig:eig_im}). The eigenimages from the traditional PCA also fail to preserve the symmetries (see e.g. 6-th eigenimage in Figure~\ref{fig:eig_im}) that are present in our eigenimages, since they do not utilize the steerable basis and rotation-augmented images.

\begin{figure}[htb]
\centering
\includegraphics[width=.75\textwidth]{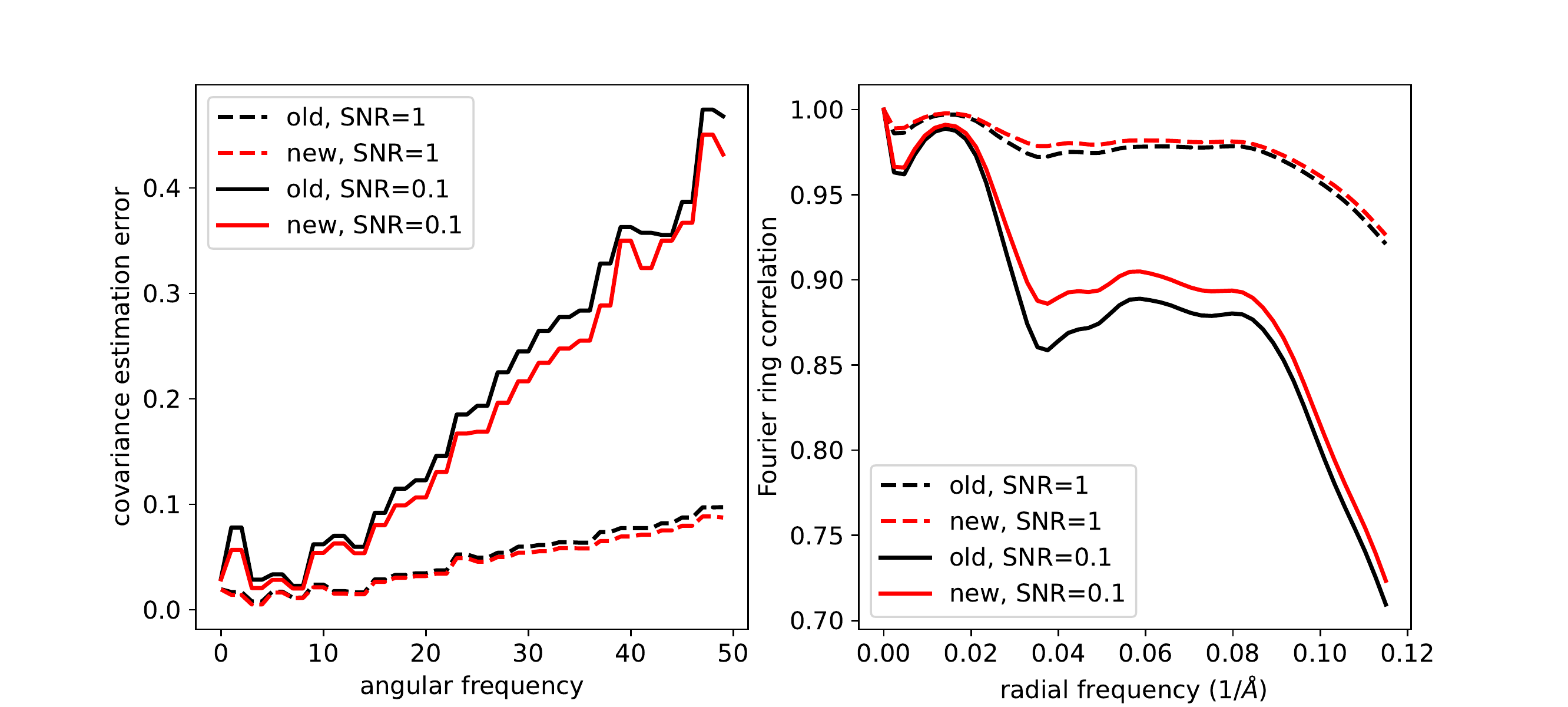}
\caption{Relative estimation error of the covariance matrix (left) and the Fourier ring correlation between the denoised and clean images (right)}\label{fig:synthetic_err}
\end{figure}

Figure~\ref{fig:synthetic_err} shows the quality of covariance estimation and image denoising when $L=128$ and using $M=100$ defocus groups. The quality of the covariance estimation is measured by the relative error in each angular frequency, which is defined as
\begin{equation}
\text{err}_n = \|C_n-\tilde{C}_n\|_F/\|C_n\|_F,
\end{equation} 
where $C_n$ and $\tilde{C_n}$ are respectively the $n$-th diagonal blocks of the clean and estimated covariance matrix, corresponding to all indices of the form $(nk,nk')$. The clean covariance matrix was approximated by the sample covariance matrix of  $10^6$ clean projection images. The performance of image denoising is measured by the Fourier ring correlation (FRC) between the clean and denoised images. Namely, for the $i$-th pair of clean and denoised images $I_i^c$ and $I_i^d$, we first compute their Fourier coefficient vectors $\boldsymbol f_{i,r}^c, \boldsymbol f_{i,r}^d$ at radial frequency $r$ by the nonuniform FFT \cite{Dutt1993,Greengard2004,lee2005type}, where $1\leq i\leq N$ and $0\leq r\leq r_{\max}$. We then compute their averaged correlation for each $r$:
$$
\text{FRC}(r) = \frac{1}{N}\sum_{i=1}^N\frac{\left\langle\boldsymbol f_{i,r}^{c}\,,  \boldsymbol f_{i,r}^d \right\rangle}{\|\boldsymbol f_{i,r}^c\|\|\boldsymbol f_{i,r}^d\|},
$$
where $\langle \cdot \,, \cdot \rangle$ denotes the inner product between two complex vectors. The FRC is a real-valued quantity due to a symmetry property that arises since the images $I_i^c, I_i^d$ are real-valued.

In addition to the speedups of Figure~\ref{fig:synth_cov}, the left panel of Figure~\ref{fig:synthetic_err} demonstrates a slight increase in the estimation quality of our proposed algorithm, compared to the previous approach. This is possibly caused by improved accuracy in the Fourier-Bessel basis approximation as well as improved accuracy by using the closed-form expression \eqref{eq:cov} compared to the approximate conjugate gradient step of previous approaches. Similarly, on the right panel, the Fourier ring correlation
between the denoised and clean images shows a slight performance increase. Figure~\ref{fig:synthetic_image} shows sample denoised images for different values of the signal-to-noise ratio (SNR) where $L=128$ and $M=100$. As a comparison, we show images denoised using the approach of this paper and the CWF method \cite{CWF}.

\section{Experimental data results}\label{sec:comp_real}
We conclude by using our method on two experimental datasets obtained from the Electron Microscopy Public Image Archive \cite{empiar}, namely EMPIAR-10028 \cite{10028}, and EMPIAR-10081 \cite{8511}. EMPIAR-10028 is a dataset of the Plasmodium falciparum 80S ribosome bound to the anti-protozoan drug emetine whose
3-D reconstruction is available in the EM data bank as EMD-2660 \cite{10028}. The dataset contains 105247 motion corrected and picked particle images, from 1081 defocus groups, of size $360\times 360$ with 1.34 \AA~pixel
size. EMPIAR-10081 is a dataset of the human HCN1 hyperpolarization-activated cyclic nucleotide-gated ion channel, whose
3-D reconstruction can be found in the EM data bank as EMD-8511 \cite{8511}. The dataset contains 55870 motion corrected and picked particle images, from 53384 defocus groups, of size $256\times 256$ with 1.3 \AA~pixel
size.

\begin{table}[t]
    \caption{Timing comparison in seconds for EMPIAR-10028 (top) and EMPIAR-10081 (bottom). Note that for EMPIAR-10081, when $L = 256$, the old-CWF method \cite{CWF} encounters memory issues and cannot be run until completion. $T_{\text{ffb}}$ denotes the time required to expand all images in the Fourier-Bessel basis (step 1 in Alg.~\ref{alg:cov}), $T_{\text{ctf}}$ the time to compute a matrix representation of the application of the point spread function (step 2 in Alg.~\ref{alg:cov}), $T_{\text{cov}}$ the time to estimate the covariance matrix (steps 3 -- 5 in Alg.~\ref{alg:cov}), $T_{\text{denoise}}$ the time to denoise the number of images indicated in the main text ($2014$ images for EMPIAR-10028 and $502$ images for EMPIAR-10081), and $T_{\text{total}}$ the total computational time.}
    \label{tab:10028}
        \centering
        \begin{tabular}{cc c c c c}
    \textbf{EMPIAR-10028} \\
    \hline
    Methods & $T_{\text{ffb}}$ & $T_{\text{ctf}}$ & $T_{\text{cov}}$ & $T_{\text{denoise}}$ & $T_{\text{total}}$\\\hline 
       old-CWF, $L=360$ & 1415  & 598 & 27550 & 201 & 29764 \\ 
       fast-CWF, $L=360$  & 768 & 5 & 2220 & 95 & 3088
       \\
       \hline
      \hline
    \end{tabular}

\begin{tabular}{cc c c c c}

    \textbf{EMPIAR-10081} \\
    \hline
    Methods & $T_{\text{ffb}}$ & $T_{\text{ctf}}$ & $T_{\text{cov}}$ & $T_{\text{denoise}}$ & $T_{\text{total}}$\\\hline 

old-CWF, $L=128$  & 47 & 3369 & 46318  & 45 & 49779\\ 
fast-CWF, $L=128$  & 35 & 8 & 93 & 27 & 163 \\\hline 
old-CWF, $L=256$  & 363 & NA& NA& NA & NA\\ 
fast-CWF, $L=256$  & 169 & 34 & 1007 & 163 & 1373
       \\
       \hline

    \end{tabular}

\end{table}

\begin{figure}[t]
\centering
\includegraphics[width=0.8\textwidth]{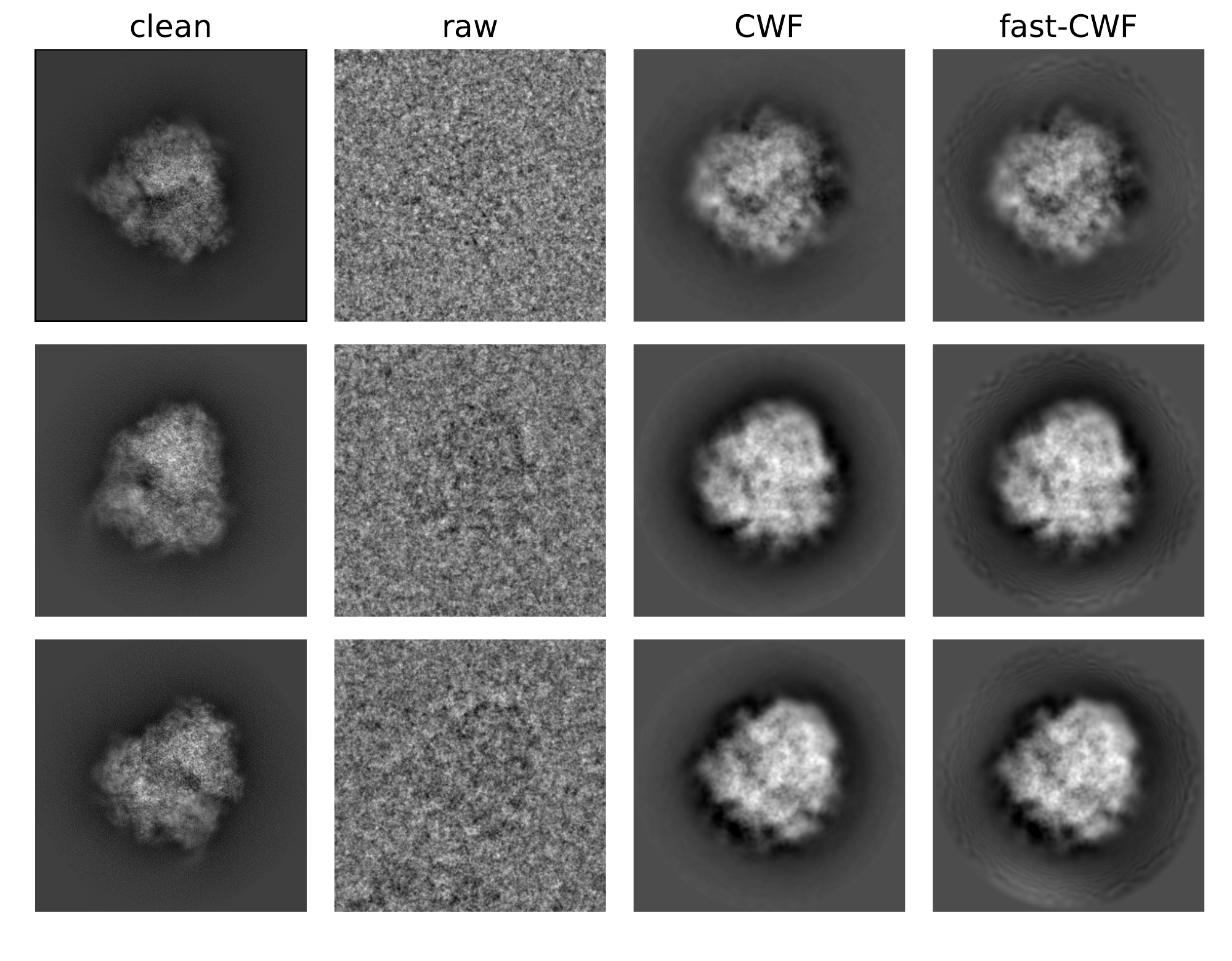}
\caption{Denoised images (EMPIAR-10028). The method used $N=105247$  images, from $M=1081$ defocus groups, of size $360\times 360$. The clean images are obtained by aligning $1000$ clean projection images (from uniformly distributed viewing directions) with phase-flipped raw images }\label{fig:denoised_10028}
\end{figure}

\begin{figure}[htb]
\centering
\includegraphics[width=0.7\textwidth]{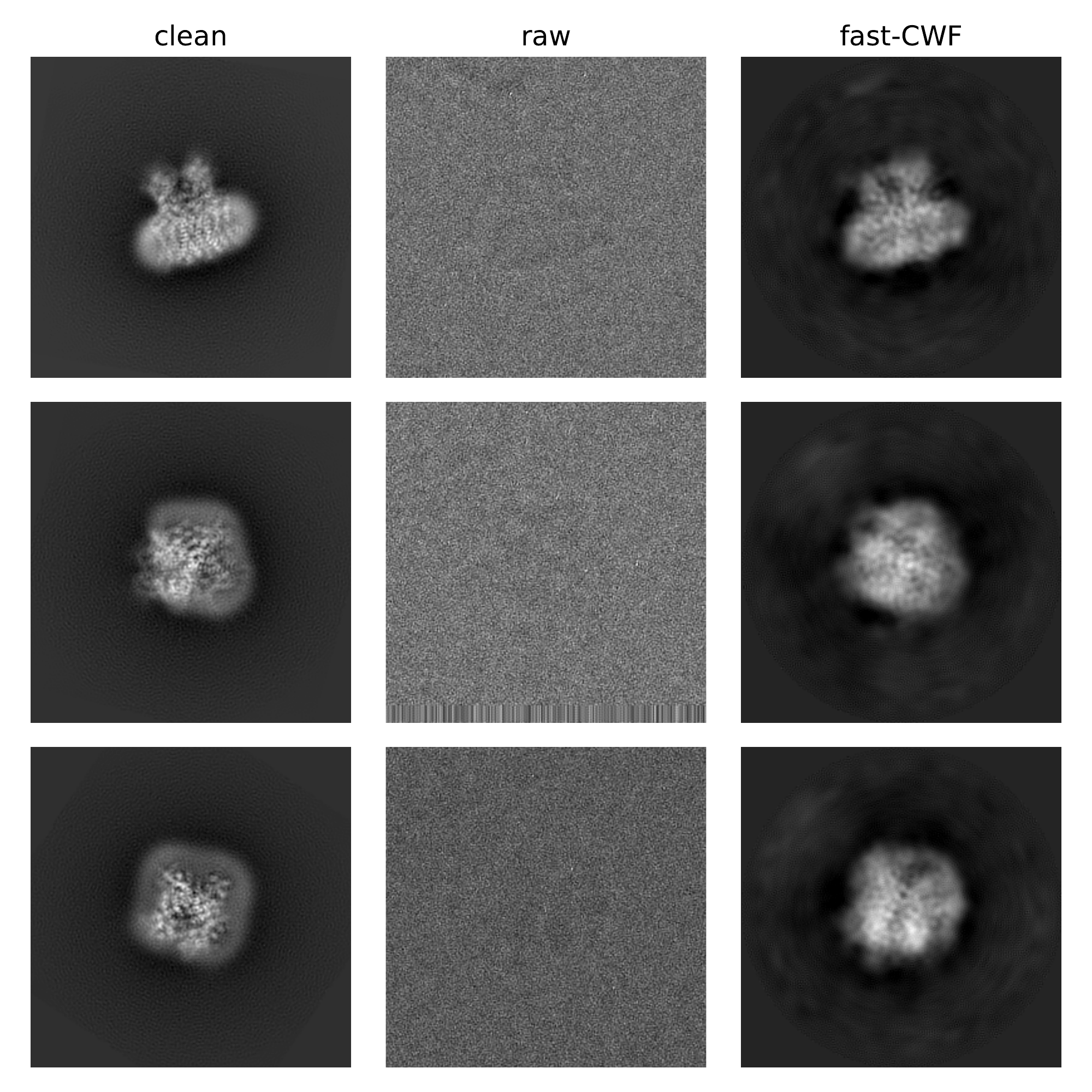}
\caption{Denoised images (EMPIAR-10081). The method used $N=55870$  images, from $M=53384$ defocus groups, of size $256\times 256$.  The clean images are obtained by aligning $1000$ clean projection images (from uniformly distributed viewing directions) with phase-flipped raw images}\label{fig:denoised_10081}
\end{figure}

Computational times are shown in Table~\ref{tab:10028}, showing a speedup of more than two orders of magnitude for the datasets with the largest number of distinct CTFs. Note that regular CWF on EMPIAR-10081 encounters memory issues and cannot be run to completion, whereas our fast method runs seamlessly. In order to obtain a comparison, we therefore additionally downsample these images to $L=128$ where the original CWF can successfully run. On EMPIAR-10028, we used all images for covariance estimation, and denoised the 2014 images from $0$-th, $50$-th, $100$-th, ..., $1000$-th defocus groups. On EMPIAR-10081, we used all images for covariance estimation, and denoised the 502 images from $0$-th, $100$-th, $200$-th, ..., $50000$-th defocus groups. For both old and new methods, the covariance matrices are further refined to correct image contrast variations \cite{shi2022ab}. Sample visualization results are shown in Figures~\ref{fig:denoised_10028} and \ref{fig:denoised_10081}.

\section{Discussion} \label{discuss}
Covariance estimation and PCA of cryo-EM images are key ingredients in many classic cryo-EM methods including multivariate statistical analysis \cite{vanHeel,van1984multivariate,van2009multivariate,van2016multivariate} and Kam's method for ab-initio modeling \cite{kam1980reconstruction}.
We propose a fast method to estimate the covariance matrix of noisy cryo-EM images, and then illustrate its application to simultaneously correct for the CTFs and denoise the images. The approach relies on recent improvements to algorithms for expanding images in the Fourier-Bessel basis \cite{FLE}, and has time complexity $O(NL^3+L^4)$ which is independent of the number of defocus groups. Our new approach is both significantly faster and more memory-efficient compared to the previous CWF method \cite{CWF} and we apply our method to large experimental datasets with many distinct CTFs with speedups by factors up to more than two orders of magnitude. 
Our approach could potentially be extended to higher-dimensional data and to the setting where images are distorted by CTFs which are not exactly radial, using either analytical correction terms or iterative numerical steps.

\bibliographystyle{abbrv}
\bibliography{cryo}

\end{document}